\documentclass[12pt,leqno]{article}
\usepackage{graphicx,psfrag,amssymb,amsmath,amsthm}

\newtheorem{theorem}{Theorem}[section]

\newtheorem{proposition}[theorem]{Proposition}

\def\beq{\begin{equation}}
\def\eeq{\end{equation}}

\def\cD{\mathcal{D}}
\def\cK{\mathcal{K}}
\def\cT{\mathcal{T}}
\def\dcD{{\partial\cD}}
\def\Tor{\text{\rm Tor}}
\def\IR{{\mathbb R}}
\def\la{\langle}
\def\ra{\rangle}

\numberwithin{equation}{section}

\begin{document}

\title{Flow-invariant hypersurfaces in semi-dispersing billiards}

\author{N. Chernov$^1$ and N. Sim\'anyi$^1$}

\date{\today}

\maketitle

\footnotetext[1]{Department of Mathematics, University of Alabama at
Birmingham, Birmingham, AL 35294. Phone 205-934-2154; Fax
205-934-9025; Email:$\ $ chernov@math.uab.edu and
simanyi@math.uab.edu.}

\begin{abstract}
This work results from our attempts to solve Boltzmann-Sinai's
hypothesis about the ergodicity of hard ball gases. A crucial
element in the studies of the dynamics of hard balls is the analysis
of special hypersurfaces in the phase space consisting of degenerate
trajectories (which lack complete hyperbolicity). We prove that if a
flow-invariant hypersurface $J$ in the phase space of a
semi-dispersing billiard has a negative Lyapunov function, then the
volume of the forward image of $J$ grows at least linearly in time.
Our proof is independent of the solution of the Boltzmann-Sinai
hypothesis, and we provide a complete and self-contained argument
here.
\end{abstract}


\section{Introduction}
\label{secI}


The ergodic hypothesis for gases of hard balls was first put forward
(in rather vague terms) by L.~Boltzmann back in the 1880's, and then
formalized by Ya.~G.~Sinai in 1963 \cite{Si63}. It states that the
gas of $N \geq 2$ identical hard balls (of small radius) on a torus
$\Tor^d$, $d \geq 2$, is ergodic, after certain necessary
reductions. The latter mean that one fixes the total energy, sets
the total momentum to zero, and restricts the center of mass to a
certain discrete lattice within the torus. The assumption of a small
radius is necessary to have the configuration space connected.

The motion of $N$ hard balls in a $d$-dimensional torus naturally
reduces to that of a single billiard particle in a
$(dN-d)$-dimensional domain $\cD$ with specular reflections off its
boundary $\dcD$. The domain $\cD$, which is the configuration space
of the system, is obtained by the removal of $N(N-1)/2$ cylinders
from the torus $\Tor^{dN-d}$ (each cylinder corresponds to a pair of
colliding balls). Thus the boundary $\dcD$ is convex (as seen from
inside the domain $\cD$), but not strictly convex, which makes the
billiard in $\cD$ semi-dispersing. This is true, unless $N=2$, in
which case we have a single `cylinder' that actually becomes a
sphere, thus $\dcD$ is strictly convex, which makes the
corresponding billiard dispersing.

The system of $N=2$ balls (disks) in dimension $d=2$ was thoroughly
investigated by Sinai \cite{Si70}: he proved its hyperbolicity,
ergodicity, and K-mixing. In this seminal paper Sinai also developed
a general theory of planar dispersing billiards. Then Sinai and
Chernov \cite{SC87} extended these results to systems of $N=2$ balls
in any dimension $d>2$, as well as to other multidimensional
dispersing billiards (there was a notable oversight in \cite{SC87}
that was corrected later in \cite{BCST02}).

The dynamics of hard balls (and related billiards) is not entirely
smooth -- singularities occur at `grazing' collisions of two balls
and at multiple collisions of $\geq 3$ balls (if they hit each other
simultaneously). These events induce hypersurfaces in the phase
space consisting of singular trajectories; they are called
singularity manifolds and cause major troubles in the analysis. In
addition, there is an issue of uniform versus nonuniform
hyperbolicity. Dispersing billiards are characterized by strong
hyperbolicity with uniform expansion and contraction rates. On the
contrary, semi-dispersing billiards are only non-uniformly
hyperbolic, so that their expansion and contraction rates may be
arbitrarily slow. This constitutes a crucial difference between
systems of $N=2$ (dispersing) and $N>2$ (semidispersing) cases, the
latter are much harder to deal with.

The first major result for systems of $N>2$ balls was the theorem on
`local' ergodicity by Sinai and Chernov \cite{SC87} (a plan of its
proof was earlier laid out by Sinai \cite{Si79}); loosely speaking
it states that for any hyperbolic phase point $x$ there is an open
neighborhood $U(x)$ that belongs to one ergodic component (mod 0).
That theorem was proven under an assumption later referred to as
Chernov-Sinai Ansatz \cite{KSS90,KSS91}; roughly speaking it states
that typical points on singularity manifolds must be hyperbolic.

A strategy for proving complete hyperbolicity and `global'
ergodicity for systems of $N > 2$ hard balls was first proposed by
Sinai and Chernov in an unpublished manuscript \cite{SC85} and
further developed and refined by Kr\'amli, Sim\'anyi, and Sz\'asz
\cite{KSS90,KSS91,KSS92}. It is based on the following observations.

A phase point $x$ fails to be hyperbolic in two possible ways.
First, it may happen that along the trajectory of $x$ the system of
$N$ balls splits into two (or more) groups (clusters) of balls that
only interact within clusters. As unlikely as it seems, such an
anomaly does occur, it is called \emph{splitting}. Second, even
without splitting the balls may collide in such a `degenerate'
manner that expansion and contraction of tangent vectors only occur
in a tangent subspace of some lower dimension. The simplest example
is obtained by letting all the $N$ balls move along a closed
geodesic in the torus (with different but parallel velocity
vectors). In such examples the balls seem to `conspire' to prevent
complete hyperbolicity; in any case their positions and velocity
vectors must satisfy very stringent requirements. We call this
phenomenon \emph{degeneracy}. Thus the complete hyperbolicity can
only be destroyed by either splitting or degeneracy.

In order to prove the hyperbolicity for the system of $N$ balls, it
is enough to verify that splitting and degeneracy occur on subsets
of zero measure. To establish ergodicity, one needs more than that
-- the set of `bad' phase points (where splitting or degeneracy
occur) must not separate two ergodic components. For instance, it is
enough to show that the set of all bad phase points has topological
codimension two or more. In addition, one needs to verify the Ansatz
to utilize the local ergodic theorem. All these proofs are carried
out inductively, by assuming the ergodicity of systems of $n$ hard
balls for all $n < N$.

The above strategy was successfully employed by Kr\'amli, Sim\'anyi
and Sz\'asz who proved hyperbolicity and ergodicity for $N=3$ balls
in any dimension \cite{KSS91} and for $N=4$ balls in dimension
$d\geq 3$, see \cite{KSS92}. Then Sim\'anyi \cite{Sim92a,Sim92b}
proved hyperbolicity and ergodicity whenever $N\leq d$ (this covers
systems with an arbitrary number of balls, but only in spaces of
high enough dimension, which is a very restrictive condition). By
using the induction hypothesis mentioned above, Sim\'anyi
\cite{Sim92a,Sim92b} also found a general argument showing that in
any system of hard balls, splitting only occurred on a subset of
zero measure and topological codimension two (or more) in phase
space, thus splitting was no longer a problem. But degeneracies
remained a key problem -- in all the above papers they were handled
by a direct analysis involving various cases specific for every
system of hard balls. This `case study' became overly complicated
with every extra ball added to the system, and it was clearly
impractical for $N>4$ balls.

Further progress required novel ideas, and Sim\'anyi and Sz\'asz
\cite{SS99} employed the methods of algebraic geometry. They assumed
that the balls had arbitrary masses $m_1, \dots, m_N$ (but the same
radius $r$). Now by taking the limit $m_N \to 0$ they were able to
reduce the dynamics of $N$ balls to the motion of $N-1$ balls, thus
utilizing a natural induction on $N$. Then, algebraic methods
allowed them to effectively analyze all possible degeneracies, but
only for typical (generic) vectors of ``external'' parameters
$(m_1,\dots,m_N,r)$; the latter needed to avoid some exceptional
submanifolds of codimension one in $\IR^{N+1}$, which remained
largely unknown. This approach led to a proof \cite{SS99} of
complete hyperbolicity (but not yet ergodicity) for all $N \geq 2$
and $d\geq 2$, and for \emph{generic} $(m_1,\dots,m_N,r)$. Later
Sim\'anyi alone upgraded the hyperbolicity to ergodicity for all
systems of $N\geq 2$ hard balls and typical (generic) values of
external parameters $(m_1,\dots,m_N,r)$; this was done in two
separate papers: \cite{Sim03} covered the case $d=2$ and
\cite{Sim04} dealt with $d>2$. Both papers again use algebraic
methods, thus their results are restricted to `generic' masses. This
was not quite satisfactory as the most interesting case of equal
masses $m_1 = \cdots = m_N$ (it is exactly this case that motivated
both Boltzmann and Sinai) remained open, as there was no guarantee
that the exceptional submanifolds would avoid it.

It is clear that methods of algebraic geometry, however powerful,
have to be abandoned if we are to obtain hyperbolicity and
ergodicity for \emph{all} gases of hard balls (including the case of
equal masses). Recently Sim\'anyi \cite{Sim02} made a partial
progress in this direction: he employed purely dynamical arguments
to establish complete hyperbolicity (but not yet ergodicity) for all
systems of hard balls (including equal masses $m_1=\cdots=m_N$). In
the course of that work, Sim\'anyi showed that degeneracies occurred
on a countable union of smooth submanifolds in the phase space with
codimension \emph{at least one}. This implies, of course, that
degeneracies are restricted to a subset of zero measure, which is
enough for hyperbolicity, but ergodicity would require those
submanifolds be of codimension \emph{at least two}.

In an attempt to find a general dynamical argument that would rule
out the existence of degeneracies on submanifolds of codimension one
(hypersurfaces) we showed that if such hypersurfaces existed, then
they would necessarily grow in size (volume) either in the past or
in the future (depending whether the Lyapunov function on them is
positive or negative). This fact turned out to be essential in a
subsequent proof of ergodicity by one of us (Sim\'anyi), under
certain conditions. Here we prove the volume growth of
hypersurfaces, as a separate fact, in the context of general
semi-dispersing billiards (which of course include all systems of
hard balls).

\section{Results and proofs}

Let $\cD$ be a bounded connected domain in $\IR^d$ or in the torus
$\Tor^d$; the boundary $\dcD$ is a finite union of $C^3$ smooth
compact hypersurfaces (each, possibly, with boundary):
$$
   \dcD = \Gamma_1 \cup \cdots \cup \Gamma_r.
$$
A billiard system is generated by a pointwise particle; it moves
freely (with a constant velocity vector) at unit speed in $\cD$ and
gets specularly reflected at the boundary $\dcD$ (by the classical
rule ``the angle of incidence is equal to the angle of
reflection'').

The phase space of the billiard system is the unit tangent bundle
$\Omega = \cD \times S^{d-1}$ over $\cD$. Thus phase points are
pairs $x = (q,v)$ where $q \in \cD$ is the position and $v \in
S^{d-1}$ the unit velocity vector of the billiard particle. The
billiard dynamics generates a flow $\Phi^t \colon \Omega \to
\Omega$. It is a Hamiltonian flow that preserves its Liouville
measure $\mu$; the latter is a direct product of uniform measures on
$\cD$ and $S^{d-1}$.

At every (regular) boundary point $q \in \Gamma_i \setminus \partial
\Gamma_i$, $1 \leq i \leq r$, we denote by $\nu(q)$ the (unique)
unit normal vector to $\Gamma_i$ pointing into $\cD$. The curvature
operator (the second fundamental form) $\cK_q$ is a self-adjoint
linear transformation acting on the tangent space $\cT_q \Gamma_i$;
is given (to the linear order) by
$$
   \nu(q+\delta q) = \nu(q) + \cK_q (\delta q).
$$
The billiard in $\cD$ is said to be dispersing (semi-dispersing) if
$\cK_q$ is positive definite (resp., positive semi-definite) at
every regular point $q \in \dcD$. Geometrically, this means that
$\dcD$ is strictly convex (resp., just convex) as seen from inside
$\cD$. For a recent detailed studies of the dynamics in
semi-dispersing billiards we refer the reader to \cite{BCST03}.

Now let $J \subset \cD$ be a small smooth compact hypersurface (a
submanifold of codimension-one with boundary), which is locally
flow-invariant, i.e.\ for every $x \in J \setminus \partial J$ there
is $\varepsilon>0$ such that $\Phi^tx \in J$ for all $|t| <
\varepsilon$. For every $x = (q,v) \in J \setminus \partial J$ we
denote by $n_x=(z,w) \in \cT_x \Omega$ a non-zero normal vector to
$J$, i.e.\ such that for any tangent vector $(\delta q,\, \delta
v)\in\cT_x\Omega$ the relation $(\delta q,\, \delta v) \in \cT_x J$
is equivalent to $\la\delta q,\,z\ra + \la\delta v,\,w\ra = 0$; here
$\langle\, \cdot\, ,\, \cdot\, \rangle$ denotes the Euclidean inner
product in $\IR^d$.

Observe that $\la w, v\ra = 0$ because $\|v\|=1$ for all $x = (q,v)
\in \Omega$. Also since $J$ is flow-invariant, we have $(v,0) \in
\cT_x J$, hence $\la z,v \ra =0$. Thus $w,z \in v^\perp$, where
$v^\perp$ denotes the hyperplane in $\IR^d$ orthogonal to $v$.

The value $Q(n_x) \colon = \la z,w \ra$ is called the
\emph{(infinitesimal) Lyapunov function}, see \cite{KB94} or part
A.4 of the Appendix in \cite{Ch94}. For a detailed account of the
relationship between the Lyapunov function $Q$, the symplectic
geometry in $\Omega$, and the dynamics, see \cite{LW95}.

\medskip\noindent{\sc Remark}.
Since the normal vector $n_x=(z,w)$ to $J$ is only determined up to
a nonzero scalar multiple, the value $Q(n_x)$ is only determined up
to a positive multiple. However, in our considerations the sign of
$Q(n_x)$ will be most important, and this sign is clearly
independent of the choice of $n_x$.
\medskip

Next, for any $t \in \IR$ the image $J_t = \Phi^t (J)$ is also a
locally flow-invariant compact hypersurface in $\Omega$. Given a
point $y_0 = (q_0, v_0) \in J \setminus \partial J$, its image $y_t
= (q_t, v_t) = \Phi^t (y_0)$ is a regular point in $J_t$ (i.e.\ $y_t
\notin \partial J_t$), unless $y_t$ is a reflection point. Given a
normal vector $n_0$ to $J$ at $y_0$, we define the unique normal
vector $n_t$ to $J_t$ at $y_t$ by
$$
    n_t=(D_{y_t} \Phi^{-t})^\ast (n_0),
$$
i.e.
\beq \label{ntdef}
    \la D_{y_0}\Phi^t (\delta y), n_t\ra =
    \la \delta y, n_0 \ra \qquad
    \forall \, \delta y \in \cT_{y_0} \Omega
\eeq
It is enough to require this for all $\delta y = (\delta q, \delta
v) \in \cT_{y_0} \Omega$ such that both $\delta q$ and $\delta y$
belong in $v_0^\perp$, because $w,s \in v_0^\perp$ as well.

We develop explicit formulas for $n_t$ ($t>0$). If there is no
collisions during an interval of time $[s,t]$ on the orbit
$\{y_t\}$, then the relation between $(\delta q_s,\, \delta v_s) \in
\cT_{y_s} \Omega$ and $(\delta q_t,\, \delta v_t) = D_{y_s}
\Phi^{t-s} (\delta q_s,\, \delta v_s)$ is obviously
$$
\aligned
  \delta v_t&=\delta v_s, \\
  \delta q_t&=\delta q_s+(t-s)\delta v_s,
\endaligned
$$
from which we obtain that for all $\delta y_s = (\delta q_s, \delta
v_s) \in \cT_{y_s} \Omega$
$$
\aligned
   \la \delta y_s, n_s \ra &=
   \la \delta q_t-(t-s)\delta v_t,\, z_s\ra
   +\langle\delta v_t,\, w_s\ra\\
   &=
   \la\delta q_t,\, z_s\ra
   +\la\delta v_t,\, w_s-(t-s) z_s\ra,
\endaligned
$$
hence
\beq \label{nbetween}
   n_t=(z_t,\, w_t)=
   \bigl(z_s,\, w_s-(t-s)z_s\bigr).
\eeq
Observe that
\beq \label{Qbetween}
   Q(n_t)=Q(n_s)-(t-s)\|z_s\|^2 \leq Q(n_s).
\eeq

Next let the orbit $\{y_t\}$ collide (reflect) at a boundary point
$q \in \dcD$ at some time $t>0$. The velocity vector $v^- = v_{t-0}$
is transformed to $v^+ = v_{t+0}$ by the rule
$$
   v^+ = v^- +2\cos\varphi \,\, \nu(q),
$$
where $\varphi$ is the angle between the outgoing velocity vector
$v^+$ and the unit normal vector $\nu(q)$ to $\dcD$, i.e.\ $\cos
\varphi = \la v^+, \nu(q) \ra$.

At the collision, the flow transforms tangent vectors $(\delta
q^-,\, \delta v^-) \in \cT_{y_{t-0}} \Omega$ (we can assume that
both $\delta q^-$ and $\delta v^-$ belong in $(v^-)^\perp$, see
above), according to the standard rules
\beq \label{RRR}
\begin{split}
   \delta q^+&=R(\delta q^-),\\
   \delta v^+&=R(\delta v^-)+2\cos\varphi \,
   RV^\ast\cK V(\delta q^-),
\end{split}
\eeq
where the operator $R \colon \cT_q \cD \to \cT_q \cD$ is the
orthogonal reflection across the tangent hyperplane $\cT_q \dcD$,
$V\colon (v^-)^\perp\to\cT_q\dcD$ is the $v^-$-parallel projection
of $(v^-)^\perp$ onto $\cT_q\dcD$, $V^\ast$ is the adjoint of $V$,
i.e.\ it is the $\nu(q)$-parallel projection of $\cT_q \dcD$ onto
$(v^-)^\perp$, and $\cK$ is the curvature operator of $\dcD$ at $q$
(see above).

The formulas (\ref{RRR}) are given, e.g., in \cite[Section 1]{SC82}
or \cite[Proposition 2.3]{KSS90}. We can also rewrite them as
$$
\aligned
  \delta q^-&=R(\delta q^+), \\
  \delta v^-&=R(\delta v^+)-2\cos\phi \,
  RV_1^\ast \cK V_1(\delta q^+),
\endaligned
$$
where $V_1$ is the $v^+$-parallel projection of $(v^+)^\perp$ onto
$\cT_q\dcD$.

Now let $n^- = (z^-, w^-)$ denote the normal vector to $J_t$ at
$y_{t-0}$, i.e.\ right before the collision, then
$$
\aligned
   \la (\delta q^-,\, \delta v^-), n^- \ra &=
   \la R(\delta q^+), \, z^-\ra+\la R(\delta v^+)-2\cos\phi
   RV_1^\ast \cK V_1(\delta q^+),\,w^-\ra \\
  &=\la\delta q^+,\, R(z^-)-2\cos\phi V_1^\ast
  \cK V_1R(w^-)\ra + \la\delta v^+,\, R(w^-)\ra.
\endaligned
$$
Thus the normal vector $n^+ = (w^+, z^+)$ to $J_t$ at $y_{t+0}$
(after the collision) is
$$
   n^+=\bigl(R(z^-)-2\cos\phi V_1^\ast \cK V_1R(w^-),\,
   R(w^-)\bigr).
$$
Observe that $\|w^+\| = \|w^-\|$ and
$$
\aligned
   Q(n^+)&=Q(n^-)-2\cos\phi\,\la V_1^\ast \cK V_1R(w^-),\, R(w^-)\ra \\
   &=Q(n^-)-2\cos\phi\,\la \cK V_1R(w^-),\, V_1R(w^-)\ra \leq
   Q(n^-),
\endaligned
$$
because $\cK \geq 0$ for semi-dispersing billiards.

We summarize our results:

\begin{proposition}
The Lyapunov function $Q(n_t)$ is a monotone non-increasing function
of $t$. The norm $\|w_t\|$ is a continuous function of $t$. The norm
$\|z_t\|$ is a piece-wise constant function of $t$ changing only at
collisions.
\end{proposition}

Our normal vector $n_t$ is related to the rate of expansion of the
hypersurface $J$ by the flow $\Phi^t$:

\begin{proposition} \label{PrVolume}
The flow $\Phi^t$ expands the volume of the hypersurface $J$ at a
point $y_0 \in J$ by a factor $\|n_t\| / \|n_0\|$.
\end{proposition}

\begin{proof}
According to (\ref{ntdef}), for every $\delta y_t \in \cT_{y_t}
\Omega$ we have
$$
   \la \delta y_t,\,n_t\ra =
   \left\la D_{y_t} \Phi^{-t}(\delta y_t),\, n_0\right\ra.
$$
Substituting $\delta y_t = n_t$ gives
$$
  \frac{\|n_t\|}{\|n_0\|}=\frac{\left\la D_{y_t}
  \Phi^{-t}(n_t),\, n_0\right \ra}{\|n_t\|\, \|n_0\|}.
$$
A simple geometric inspection shows that the right-hand-side is
exactly the factor of linear expansion by $\Phi^{-t}$ between the
points $y_t$ and $y_0$ in the direction transversal to $J_t$, and
since $\Phi^t$ preserves volume in $\Omega$, it is the factor of
volume expansion of the hypersurface $J$ by $\Phi^t$ between the
points $y_0$ and $y_t$. The proposition is proved.
\end{proof}

From now on we assume that $Q (n_0)<0$. This assumption is motivated
by the following considerations (which will not affect our further
arguments, though). If there is a flow-invariant compact
hypersurface $J \subset \Omega$, it might divide the phase space
$\Omega$ into two (or more) non-interacting open subsets, thus
preventing the ergodicity of the flow. One wants to show that such
hypersurfaces cannot exist. There is a separate argument
\cite[Remark 7.9]{Sim03} showing that the Lyapunov function cannot
stay identically zero on the entire trajectory of any non-splitting
phase point, thus $Q<0$ or $Q>0$ on parts of $J$. We will show that
if $Q<0$, then the volume of $J$ is expanded under the flow, at
least linearly in time, hence $J$ cannot be compact. Similarly, if
$Q>0$, then the volume of $J$ is expanded in the past.

\begin{proposition}
If $Q(n_0)<0$, then $Q(n_t)$ is a monotone decreasing function of
$t$, and the norm $\|w_t\|$ of the `velocity' component of the
normal vector $n_t=(z_t,\,w_t)$ is a monotone increasing function of
$t$.
\end{proposition}

\begin{proof}
We have $Q(n_s) = \la z_s, w_s \ra \leq Q(n_0) <0$ for all $s>0$.
Hence $\|z_s\| \neq 0$, and so $Q(n_t) < Q(n_s)$ for any
collision-free time interval $(s,t)$ due to (\ref{Qbetween}). Also,
by (\ref{nbetween})
$$
  \|w_t\|^2 = \|w_s\|^2 - 2(t-s)\, Q(n_s)
  +(t-s)^2\|z_s\|^2 > \|w_s\|^2,
$$
and $\|w_t\|$ does not change at collisions.
\end{proof}

\begin{proposition}
If $Q(n_0)<0$, then the function $\|w_t\|/|Q(n_t)|$ is monotone
decreasing in $t$.
\end{proposition}

\begin{proof}
For any $s\geq 0$ we have
$$
   w_s = \frac{Q(n_s)}{\|z_s\|^2}\, z_s + w_s^\perp,
$$
where $w_s^\perp$ is a vector orthogonal to $z_s$, hence
$$
   \frac{\|w_s\|}{|Q(n_s)|} =
   \biggl[\frac{1}{\|z_s\|^2} + \frac{\|w_s^\perp\|^2}
   {|Q(n_s)|^2}\biggr]^{1/2}
$$
Due to (\ref{nbetween}), for any collision-free time interval
$(s,t)$ we have $z_t = z_s$ and $w_t^\perp = w_s^\perp$, while
$|Q(n_s)| < |Q(n_t)|$. At collisions, $\|w_t\|$ stays constant,
while $|Q(n_t)|$ instantaneously increases, thus the required
property holds.
\end{proof}

\begin{proposition} \label{PrLinear}
If $Q(n_0) <0$, then for all $t\geq 0$
$$
  \|w_t\|\geq \|w_0\|+\frac{|Q(n_0)|\, t}{\|w_0\|}.
$$
\end{proposition}

\begin{proof}
For any $t>0$, which is not a collision time,
$$
  \frac{d}{dt}\|w_t\|^2=\frac{d}{dt}\la w_t,\,w_t\ra=-2
  \la z_t,\,w_t\ra=2|Q(n_t)|\geq\frac{2\, |Q(n_0)|\,\|w_t\|}{\|w_0\|},
$$
where we used (\ref{nbetween}) and the previous proposition, hence
$$
  \frac{d}{dt}\|w_t\|\geq\frac{|Q(n_0)|}{\|w_0\|}.
$$
Since $\|w_t\|$ does not change at collisions, the proposition
follows.
\end{proof}

Next we assume that $Q(t) < 0$ on all of the manifold $J$; more
specifically
\beq \label{c0}
  Q(n_0)=\la z_0,\,w_0\ra \leq -c_0<0
\eeq
holds true uniformly across $J$ for a unit normal vector field
$n_0(y)$ ($y\in J$) with a fixed constant $c_0>0$.

\begin{theorem}
If (\ref{c0}) holds, then $\Phi^t$ expands the volume of the
hypersurface $J$ at any point $y \in J$ by a factor $\lambda_t \geq
1+c_0t$ for all $t \geq 1/c_0$.
\end{theorem}

\begin{proof}
Combining Propositions~\ref{PrVolume} and \ref{PrLinear} gives
$$
 \lambda_t = \|n_t\| \geq \|w_t\| \geq \|w_0\|+\frac{c_0t}{\|w_0\|}.
$$
Since $\|w_0\| \leq \|n_0\| = 1$, this lower bound is at least
$1+c_0t$ for $t\ge 1/c_0$.
\end{proof}

\textbf{Acknowledgement}. N.~C.\ was partially supported by NSF
grant DMS-0354775. N.~S.\ was partially supported by NSF grant
DMS-0457168.

\end{document}